\documentclass[12pt,reqno]{amsart}

\usepackage{graphicx}
\usepackage[draft]{hyperref}
\usepackage{amsmath,amsopn,amssymb,amsfonts,stmaryrd}
\usepackage{verbatim}
\usepackage{amsthm}
\usepackage{mathtools}
\usepackage{color}
\usepackage{enumitem}
\usepackage[framemethod=TikZ]{mdframed}
\usepackage{bbm}
\usepackage{mathrsfs}
\usepackage{booktabs}
\usepackage{caption}
\usepackage{bm}

\usepackage{tensor}
\usepackage{xcolor}
\usepackage{bbm}
\usepackage{cleveref}

\makeatletter
\@namedef{subjclassname@2020}{%
 \textup{2020} Mathematics Subject Classification}
\makeatother

\newcommand{\Q}{\mathbb{Q}}

\newcommand{\SU}{\mathfrak{su}}

\newcommand{\su}[1]{{\mathfrak{su}{#1}}}
\newcommand{\so}[1]{{\mathfrak{so}{#1}}}

\newcommand{\N}{\mathbb N}
\newcommand{\C}{\mathbb C}

\setlist[enumerate]{leftmargin=*,label=\rm{(\arabic*)}}

\theoremstyle{plain}
\newtheorem{thm}{Theorem}[section]

\theoremstyle{definition}

\numberwithin{equation}{section}

\theoremstyle{theorem}

\newcommand{\flo}[1]{\lfloor #1\rfloor}
\newcommand{\Flo}[1]{\left\lfloor #1\right\rfloor}
\newcommand{\pa}[2]{\left(\frac{#1}{#2}\right)}
\newcommand{\sm}{\setminus}

\def\p{\varrho}

\def\w{\omega}

\def\z{\zeta}

\def\g{\gamma}

\def\p{\varrho}

\def\w{\omega}

\def\z{\zeta}

\def\g{\gamma}

\def\GG{\Gamma}

\newcommand{\re}{{\rm Re}}

\setlist[itemize]{leftmargin=*}

\makeatletter
\newcommand{\vast}{\bBigg@{2}}
\newcommand{\Vast}{\bBigg@{5}}
\makeatother

\newcommand{\Pmod}[1]{\ \, ( \mathrm{mod} \, #1 )}

\newcommand{\Res}{\operatorname{Res}}

\definecolor{Green}{rgb}{0,0.4,0}

\allowdisplaybreaks

\title{On the number of irreducible representations of $\su(3)$}

\author{Walter Bridges}
\author{Kathrin Bringmann}
\author{Johann Franke}
\address{University of Cologne, Department of Mathematics and Computer Science, Weyertal 86-90, 50931 Cologne, Germany}
\email{wbridges@uni-koeln.de}
\email{kbringma@math.uni-koeln.de}
\email{jfrank12@uni-koeln.de}



\begin{document}
\maketitle

\begin{abstract}
	In this note, we use a variant of the hyperbola method to prove an asymptotic expansion for the summatory function of the number of irreducible $\su(3)$-representations of dimension $n$. This is a natural companion result to work of Romik, who proved an asymptotic formula for the number of unrestricted $\su(3)$-representations of dimension $n$.
\end{abstract}

\section{Introduction and statement of results}

The irreducible representations of the Lie algebra $\su(3)$ are a family of representations $W_{j,k}$ of dimension $\frac{jk(j+k)}{2}$ for $j,k\in\N_0$ (see \cite[Theorem 6.27]{Hall}). Let $r(n)$ denote the number of $n$-dimensional representations of $\su(3)$. Then
\begin{equation}\label{E:rngenfn}
	\sum_{n \geq 0} r(n)q^n = \prod_{j,k\geq 1} \frac{1}{1-q^{\frac{jk(j+k)}{2}}}.
\end{equation}
Romik recently proved the following asymptotic formula for $r(n)$ by studying (a renormalization of) the Witten zeta function for $\mathrm{SU}(3)$; that is, the meromorphic continuation of the series
\begin{align*}
	\z_{\su(3)}(s) := \sum_{j,k \geq 1} \left( \frac{jk(j+k)}{2}\right)^{-s} \qquad \left(\re(s) > \tfrac23\right).
\end{align*}

\begin{thm}[\cite{Ro}, Theorem 1.1]\label{T:Romik}
	As $n\to\infty$ , we have, for certain constants $A_1,A_2,A_3,A_4,K>0$
	$$
		r(n) \sim \frac{K}{n^{\frac{3}{5}}}\exp\left(A_1n^{\frac{2}{5}}-A_2n^{\frac{3}{10}}-A_3n^{\frac{1}{5}}-A_4n^{\frac{1}{10}} \right).
	$$
\end{thm}

Romik stated that Theorem \ref{T:Romik} is an analogue of the Hardy--Ramanujan asymptotic formula for $p(n)$, the number of integer partitions of $n$, because the corresponding generating function for $\su(2)$-representations coincides with
\begin{equation}\label{E:pngenfn}
	\sum_{n\ge0} p(n)q^n = \prod_{n\ge1} \frac{1}{1-q^n}.
\end{equation}
The doubly indexed product \eqref{E:rngenfn} has much more complicated analytic behavior compared to the modular infinite product \eqref{E:pngenfn}. Two of the authors \cite{BringmannFranke} subsequently obtained an asymptotic series for $r(n)$ which was then generalized by the authors and Brindle to more general product generating functions \cite{BBBF}, including for example representations of $\so(5)$.

In the present paper, we turn our attention to the number of irreducible $\SU(3)$-representations of dimension $n$, i.e.,
\begin{align*}
	\p(n) = \sum_{\substack{j,k \geq 1 \\ \frac{jk(j+k)}{2} = n}} 1.
\end{align*}
Of course, this is a highly oscillatory function and is often 0, but we may still study the average, $\sum_{1 \leq n\leq x} \varrho(n),$ as $x \to \infty$. Note the similarity of $\varrho(n)$ to the divisor function,
$$
	d(n):=\sum_{\substack{j,k\geq 1 \\ jk=n}}1.
$$
Dirichlet's hyperbola method yields the first two terms in the expansion of the average,
$$
	\sum_{1 \leq n \leq x} d(n)=x\log(x)+(2\gamma -1)x + O\left(\sqrt{x}\right),
$$
where $\g$ is the Euler--Mascheroni constant (see for example \cite[Theorem 3.3]{Apostol}). The still open Dirichlet divisor problem concerns improving the error term from $O(\sqrt x)$ to the conjectured $O(x^\frac14)$; for an overview, see \cite{Berndt}.

We show here that a variant of the hyperbola method yields the following asymptotic expansion for the summatory function of $\varrho(n)$.

\begin{thm}\label{T:Main}
	We have, as $x\to\infty$,
	\[
		\sum_{1\le n\le x} \p(n) = \frac{2^\frac23\sqrt3\GG\pa13^3}{4\pi}x^\frac23 + 2^\frac32\z\pa12\sqrt{x} + O\left(x^\frac13\right).
	\]
\end{thm}

We prove \Cref{T:Main} in \Cref{S:Proof}, and we conclude this section with the following questions.
\begin{enumerate}
	\item Can one improve the error term in \Cref{T:Main}, perhaps by a deeper study of the Witten zeta function, $\z_{\su(3)}(s)$? It would be reasonable to consider deeper techniques that have been brought to bear on the Dirichlet divisor problem (the Selberg--Delange method \cite{Tenenbaum}, Voronoi summation \cite{IK}, to name a few).
	
	\item Can this variant of the hyperbola method (or any other technique) be used to yield asymptotic series for generic sums
	$$
		\sum_{\substack{m,n\geq 1 \\ p(m,n)\leq x}} 1,
	$$
	where $p(x,y)$ is a homogeneous polynomial in $\Q[x,y]$ taking integer values? For example, the case $p(m,n)=\frac{mn(m+n)(m+2n)}{6}$ corresponds to representations of $\so(5)$.
\end{enumerate}

\section*{Acknowledgements}

We thank Dan Romik for sharing this problem with us and for providing helpful feedback. The authors have received funding from the European Research Council (ERC) under the European Union's Horizon 2020 research and innovation programme (grant agreement No. 101001179).

\section{Proof of Theorem \ref{T:Main}}\label{S:Proof}

In this section we prove \Cref{T:Main}.
\begin{proof}[Proof of \Cref{T:Main}]
	We note that the asymptotic main term in \Cref{T:Main} may be obtained by analytic properties of $\z_{\su(3)}(s)$ along with a standard Tauberian theorem.  In particular, Theorem 1.2 (3) of \cite{Ro} implies\footnote{Romik defined $\w(s):=2^{-s}\z_{\su(3)}(s)$.} that
\[
	\Res_{s=\frac23} \z_{\su(3)}(s) = \Res_{s=\frac23} 2^s\w(s) = \frac{2^\frac23\GG\pa13^3}{2\pi\sqrt3}.
\]
By the same theorem, $\z_{\su(3)}$ has a meromorphic continuation to $\C\sm(\frac23\cup(\frac12-\N_0))$. It follows from the Ikehara--Wiener Tauberian theorem (see, e.g. \cite[Ch. 3, ex. 3.3.6]{Murty})
\begin{equation}\label{E:mainterm}
	\sum_{1\le n\le x} \p(n) \sim \frac32\Res_{s=\frac23} \z_{\su(3)}(s)x^\frac23 = \frac{2^\frac23\sqrt3\GG\pa13^3}{4\pi}x^\frac23,\qquad x \to \infty.
\end{equation}
	
For the next term in the asymptotic expansion, we first write
\[
	\sum_{1\le n\le x} \p(n) = \sum_{1\le N\le x} \sum_{\substack{m,n\ge1\\\frac{mn(m+n)}{2}=N}} 1 = \sum_{\substack{m,n\ge1\\mn(m+n)\le2x\\mn(m+n)\equiv0\Pmod2}} 1.
\]
Now $mn(m+n)\equiv0\Pmod2$ is automatically satisfied, and we see that the above sum counts lattice points in the $(m,n)$-plane between $m=1$, $n=1$ and the curve $n=\frac{-m^2+\sqrt{m^4+8mx}}{2m}$ (the positive solution to the quadratic equation $mn(m+n)=2x$). In usual hyperbola-method fashion, we add up the lattice points for each $1\le m\le x^\frac13$ along the vertical lines from $n=1,\dots,\flo{\frac{-m^2+\sqrt{m^4+8mx}}{2m}}$. Then we do the same for $1\le n\le x^\frac13$ along the horizontal lines from $m=1,\dots,\flo{\frac{-n^2+\sqrt{n^4+8nx}}{2n}}$. By symmetry these are the same. Then we subtract the points counted twice, namely those in the square with side length $x^\frac13$. The result is
\begin{align}\nonumber
	\sum_{1\le n\le x} \p(n) &= 2\sum_{1\le n\le x^\frac13} \sum_{1\le m\le\frac{-n^2+\sqrt{n^4+8nx}}{2n}} 1 - \Flo{x^\frac13}^2\\
	\nonumber
	&= 2\sum_{1\le n\le x^\frac13} \Flo{\frac{-n^2+\sqrt{n^4+8nx}}{2n}} - x^\frac23 + O\left(x^\frac13\right)\\
	\nonumber
	&= \sum_{1\le n\le x^\frac13} \left(-n+\sqrt{n^2+\frac{8x}{n}}\right) - x^\frac23 + O\left(x^\frac13\right)\\
	\label{E:hyperbolasymplified}
	&= -\frac32x^\frac23 + \sum_{1\le n\le x^\frac13} \sqrt{n^2+\frac{8x}{n}} + O\left(x^\frac13\right).
\end{align}
Thus, we only need to approximate the remaining sum. Let $\{t\}:=t-\flo t$.  Abel partial summation \cite[Theorem 0.3, p. 4]{Tenenbaum} gives
\begin{align*}
	\sum_{1\le n\le x^\frac13} \sqrt{n^2+\frac{8x}{n}} &= \Flo{x^\frac13}\sqrt{x^\frac23+8x^\frac23} - \frac12\int_1^{x^\frac13} \frac{2t-\frac{8x}{t^2}}{\sqrt{t^2+\frac{8x}{t}}}\flo t dt\\
	&= 3x^\frac23 - \int_1^{x^\frac13} \frac{t-\frac{4x}{t^2}}{\sqrt{t^2+\frac{8x}{t}}}(t-\{t\}) dt + O\left(x^\frac13\right)\\
	&= 3x^\frac23- \int_1^{x^\frac13} \frac{t^2-\frac{4x}{t}}{\sqrt{t^2+\frac{8x}{t}}} dt + \int_1^{x^\frac13} \frac{t-\frac{4x}{t^2}}{\sqrt{t^2+\frac{8x}{t}}}\{t\} dt + O\left(x^\frac13\right).
\end{align*}
Making the change of variables $t\mapsto2x^\frac13t$, the first integral is
\[
	\int_1^{x^\frac13} \frac{t^2-\frac{4x}{t}}{\sqrt{t^2+\frac{8x}{t}}} dt = 2x^\frac23\int_\frac{1}{2x^\frac13}^\frac12 \frac{2t^\frac52-t^{-\frac12}}{\sqrt{1+t^3}} dt =: 2x^\frac23F\pa{1}{2x^\frac13},
\]
where
$$
	F(y):=\int_y^{\frac{1}{2}}\frac{2t^{\frac{5}{2}}-t^{-\frac{1}{2}}}{\sqrt{1+t^3}}dt.
$$
Noting that $\frac{1}{\sqrt{1+t^3}}=1+O(t^3)$ gives the expansion
$$
	F(y)=F(0)+2\sqrt{y}+O\left(y^{\frac{7}{2}}\right), \qquad (y \to 0).
$$
Hence,
\begin{align}\nonumber
	\sum_{1\le n\le x^\frac13} \sqrt{n^2+\tfrac{8x}{n}} &= 3x^\frac23 - 2x^\frac23\left(F(0)+\sqrt2x^{-\frac16}+O\left(x^{-\frac76}\right)\right)\\
	\nonumber
	&\hspace{6cm}+ \int_1^{x^\frac13} \tfrac{t-\frac{4x}{t^2}}{\sqrt{t^2+\frac{8x}{t}}}\{t\} dt + O\left(x^\frac13\right)\\
	\label{E:sqrtsum}
	&= (3-2F(0))x^\frac23 - 2\sqrt{2x} + \int_1^{x^\frac13} \tfrac{t-\frac{4x}{t^2}}{\sqrt{t^2+\frac{8x}{t}}}\{t\} dt + O\left(x^\frac13\right).
\end{align}
The integral with the fractional part is
\begin{align}
	\int_1^{x^\frac13} \frac{t-\frac{4x}{t^2}}{\sqrt{t^2+\frac{8x}{t}}}\{t\} dt
	&= \int_1^{x^\frac13} \left(\frac{t^\frac32\{t\}}{\sqrt{t^3+8x}} - \frac{4x\{t\}}{t^\frac32\sqrt{t^3+8x}}\right) dt. \label{E:fracpartint}
\end{align}
Now, the function $t \mapsto \frac{t^\frac32}{\sqrt{t^3+8x}}$ is increasing for $t >0$, so
\[
	\int_1^{x^\frac13} \frac{t^\frac32\{t\}}{\sqrt{t^3+8x}} dt \le \int_1^{x^\frac13} \frac{t^\frac32}{\sqrt{t^3+8x}}dt \le x^\frac13\max_{1\le t\le x^\frac13} \frac{t^\frac32}{\sqrt{t^3+8x}} = O\left(x^\frac13\right).
\]
The second integral in \eqref{E:fracpartint} is 
\[
	-4x \int_1^{x^{\frac{1}{3}}}\frac{\{t\}}{t^\frac32\sqrt{t^3+8x}}dt = -\sqrt{2x}\int_1^{x^{\frac{1}{3}}} \frac{\{t\}}{t^\frac32\sqrt{1+\frac{t^3}{8x}}}dt.
\]
Now, we have\footnote{We write $f(y)=O_{\le c}(g(y))$ for $|f(y)|\le c|g(y)|$.} $(1+y)^{-\frac12}=1+O_{\le0.5}(y)$ for $0\leq y\le\frac18$. Thus,
\begin{align*}
	\int_1^{x^\frac13} \frac{\{t\}}{t^\frac32\sqrt{1+\frac{t^3}{8x}}} dt &= \int_1^{x^\frac13} \frac{\{t\}}{t^\frac32} dt + \int_1^{x^\frac13} \frac{\{t\}}{t^\frac32}O_{\le0.5}\pa{t^3}{8x} dt\\
	&= \int_1^\infty \frac{\{t\}}{t^\frac32} dt + O\left(x^{-\frac16}\right).
\end{align*}
Now, by \cite[p. 232]{Tenenbaum} we have
$$
	\zeta\left(\frac{1}{2}\right)=-1-\frac{1}{2}\int_1^{\infty} \frac{\{t\}}{t^{\frac{3}{2}}}dt,
$$
so
\[
	\int_1^{x^{\frac{1}{3}}}\frac{\{t\}}{t^\frac32\sqrt{1+\frac{t^3}{8x}}} dt=-2-2\zeta\left(\frac{1}{2}\right)+O\left(x^{-\frac{1}{6}}\right).
\]
Thus,
\[
	-4x\int_1^{x^\frac13} \frac{\{t\}}{t^\frac32\sqrt{t^3+8x}}dt = 2\sqrt2\left(1+\z\pa12\right)\sqrt x + O\left(x^\frac13\right).
\]
Plugging into \eqref{E:sqrtsum}, we get
\[
	\sum_{1\le n\le x^\frac13} \sqrt{n^2+\frac{8x}{n}} = (3-2F(0))x^\frac23 + 2^\frac32\z\pa12\sqrt x + O\left(x^\frac13\right),
\]
which if added to \eqref{E:hyperbolasymplified} yields
\[
	\sum_{1\le n\le x} \p(n) = \left(\frac32-2F(0)\right)x^\frac23 + 2^\frac32\z\pa12\sqrt x + O\left(x^\frac13\right).
\]
Comparing with \eqref{E:mainterm}, we conclude \Cref{T:Main}.
\end{proof}
\noindent \textbf{Remark.} 
	The proof of \Cref{T:Main} implies the identity
\[
	\int_0^\frac12 \frac{2t^\frac52-t^{-\frac12}}{\sqrt{1+t^3}} dt = \frac34 - \frac{2^\frac23\sqrt3\GG\pa13^3}{8\pi}.
\]
We did not find a direct proof, but we note that the factor $\frac{\Gamma\left( \frac{1}{3}\right)^3}{\pi}$ appears in evaluations of the complete elliptic integral of the first kind \cite[Table 9.1]{BB}.

\end{document}